\documentclass[11pt,a4paper]{amsart}
\usepackage{amsmath,amssymb}
\newcounter{teo}[section]
\newcounter{def}[section]
\newcounter{prop}[section]
\newcounter{lem}[section]
\newcounter{ex}[section]
\newcounter{rem}[section]
\newcounter{cor}[section]
\numberwithin{equation}{section} \numberwithin{def}{section}
\numberwithin{prop}{section} \numberwithin{teo}{section}
\numberwithin{lem}{section} \numberwithin{ex}{section}
\numberwithin{rem}{section} \numberwithin{cor}{section}

\newenvironment{theorem}[1][\ ]
  {\newcommand{\nume}{\arabic{section}.\arabic{teo}}
   \par\bigskip\noindent\refstepcounter{teo}\textbf{{Th\'{e}or\`{e}me\ }\nume\ #1}\bgroup\it}
  {\egroup\par}
\newenvironment{definition}[1][\ ]
  {\newcommand{\nume}{\arabic{section}.\arabic{def}}
   \par\bigskip\noindent\refstepcounter{def}\textbf{{D\'{e}finition\ }\nume\ #1}\bgroup}
  {\egroup\par}
\newenvironment{proposition}[1][\ ]
  {\newcommand{\nume}{\arabic{section}.\arabic{prop}}
   \par\bigskip\noindent\refstepcounter{prop}\textbf{{Proposition\ }\nume\ #1}\bgroup\it}
  {\egroup\par}

\newenvironment{example}[1][\ ]
  {\newcommand{\nume}{\arabic{section}.\arabic{ex}}
   \par\bigskip\noindent\refstepcounter{ex}\textbf{{Exemple\ }\nume\ #1}\bgroup\it}
  {\egroup\par}
\newenvironment{remark}[1][\ ]
  {\newcommand{\nume}{\arabic{section}.\arabic{rem}}
  \par\bigskip\noindent\refstepcounter{rem}\textbf{{Remarque\ }\nume\ #1}\bgroup\it}
  {\egroup\par}
\newenvironment{corollary}[1][\ ]
  {\newcommand{\nume}{\arabic{section}.\arabic{cor}}
  \par\bigskip\noindent\refstepcounter{cor}\textbf{{Corollaire\ }\nume\ #1}\bgroup\it}
  {\egroup\par}

\begin{document}

\title[Instabilit\'{e} des cocycles d'\'{e}volution]{Instabilit\'{e} des cocycles d'\'{e}volution fortement mesurables dans des espaces de Banach }
\author{Codru\c{t}a Stoica}
\date{}
\maketitle

{\footnotesize \noindent \textbf{R\'{e}sum\'{e}.} Le but de cet
\'{e}tude est de pr\'{e}senter quelques comportements
asymptotiques concernant les cocycles d'\'{e}volution dans des
espaces de Banach, comme, par exemple, la d\'{e}croissance
exponentielle, l'instabilit\'{e}, l'instabilit\'{e} exponentielle
et l'instabilit\'{e} int\'{e}grale. Des relations entre ces
propri\'{e}t\'{e}s sont aussi d\'{e}montr\'{e}es. Deux
th\'{e}or\`{e}mes de type Datko
sont \'{e}nonc\'{e}s et d\'{e}montr\'{e}s comme r\'{e}sultats principaux. 
En m\^{e}me temps un outil d'\'{e}tude non uniforme est fourni.}

{\footnotesize \vspace{3mm} }


{\footnotesize \vspace{2mm} }

{\footnotesize \noindent \textit{Mots-clefs:} cocycle
d'\'{e}volution fortement mesurable, d\'{e}croissance
exponentielle, instabilit\'{e} exponentielle, instabilit\'{e}
int\'{e}grale}

\vspace{3mm}

{\footnotesize \noindent \textbf{Abstract.} The aim of the paper
is to present various asymptotic behaviors of skew-evolution
semiflows in Banach spaces, as exponential decay, instability,
exponential instability and integral instability. Relations
between these asymptotic properties are also given. As main
results, two Datko type theorems are proved. A unified nonuniform
approach is provided.}

{\footnotesize \vspace{2mm} }


\maketitle

\section{Introduction}

Dans l'approfondissement de l'\'{e}tude d'un cas particulier d'un
concept de cocycle d'\'{e}volution introduit dans l'ouvrage
\cite{MeSt_CJM}, concernant les \'{e}quations d'\'{e}volution,
l'approche par la th\'{e}orie des op\'{e}rateurs d'\'{e}volution a
\'{e}t\'{e} essentielle. La notion s'est montr\'{e}e plus
appropri\'{e}e dans le cas non uniforme des comportements
asymptotiques. Des propri\'{e}t\'{e}s asymptotiques pour le cas particulier de
cocycles d'\'{e}volution ont \'{e}t\'{e} pr\'{e}sent\'{e}es dans
l'article \cite{MeSt_CJM} et aussi dans \cite{MeSt_AUVT}.
L'instabilit\'{e} a \'{e}t\'{e} \'{e}tudi\'{e}e dans \cite{MeSt_ICNODEA} dans le cadre uniforme pour des op\'{e}rateurs
d'\'{e}volution. Des r\'{e}sultats concernant la stabilit\'{e} et 
l'instabilit\'{e} ont \'{e}t\'{e} li\'{e}s \`{a} l'\'{e}tude de la 
trichotomie dans \cite{MeSt_IEOT} dans le cas des op\'{e}rateurs
d'\'{e}volution.

Quelques caract\'{e}risations pour la propri\'{e}t\'{e} de
stabilit\'{e} et d'instabilit\'{e} pour les cocycles relatifs aux
semiflots ont \'{e}t\'{e} pr\'{e}sent\'{e}es dans \cite{MeSaSa}.
En g\'{e}n\'{e}ral, dans ce cas l\`{a}, aborder le probl\`{e}me du
point de vue des semigroupes d'op\'{e}rateurs a jou\'{e} un
r\^{o}le important. 

Nous envisageons de pr\'{e}senter dans cet article de nouvelles d\'{e}finitions pour les
notions de semiflot d'\'{e}volution et de cocycle d'\'{e}volution, ainsi que la 
propri\'{e}t\'{e} d'instabilit\'{e} exponentielle, trait\'{e}e
dans le cas non uniforme. Nous faisons la remarque que l'approche
est differente de la caract\'{e}risation de ce comportement
asymptotique faite dans \cite{MeSt_CJM}, o\`{u} dans la
d\'{e}finition d'un cocycle d'\'{e}volution nous avons suppos\'{e}
aussi la croissance exponentielle. Au contraire, dans l'\'{e}tude 
ci-pr\'{e}sente on utilise la d\'{e}croissance exponentielle et la
propri\'{e}t\'{e} d'instabilit\'{e} int\'{e}grale. Notre intention est
d'\'{e}tendre gr\^{a}ce aux Th\'{e}or\`{e}me
\ref{caract_uiis_1} et Th\'{e}or\`{e}me \ref{caract_uiis_2} le
r\'{e}sultat classique et bien connu contenu dans le
Th\'{e}or\`{e}me $11$ de l'ouvrage \cite{Da}, concernant le cas de
la stabilit\'{e} uniforme exponentielle. Nos d\'{e}monstrations ne
sont pas de simples g\'{e}n\'{e}ralisations parce que dans notre
\'{e}tude les cocycles d'\'{e}volution consid\'{e}r\'{e}s sont
d'un type plus g\'{e}n\'{e}ral et pas n\'{e}cessairement fortement
continus. Le cadre non uniforme est plus g\'{e}n\'{e}ral et les
r\'{e}sultats obtenus peuvent \^{e}tre appliqu\'{e}s pour des
caract\'{e}risations similaires pour le cas d'op\'{e}rateurs
d'\'{e}volution.

\section{Notations et d\'{e}finitions}

On consid\`{e}re un espace m\'{e}trique $\mathcal{X}$, un espace
de Banach $\mathcal{V}$. Soit $\mathcal{B}(\mathcal{V})$ l'espace
de tous les op\'{e}rateurs born\'{e}s,
d\'{e}finis sur $\mathcal{V}$ \`{a} valeurs dans $\mathcal{V}$. On note par $\left\Vert \cdot
\right\Vert$ la norme des vecteurs sur $\mathcal{V}$ et des
op\'{e}rateurs sur $\mathcal{B}(\mathcal{V})$.

Soit $\mathcal{T}$ l'ensemble de toutes les paires $(t,s)$ des
nombres r\'{e}els non n\'{e}gatifs tels que $t\geq s$. On
note $\mathcal{Y}=\mathcal{X}\times \mathcal{V}$ et $I$
l'op\'{e}rateur identit\'{e} sur $\mathcal{V}$.

\begin{definition}\label{semiflow}\rm
Une application $\varphi:\mathbb{R}_{+}\times
\mathcal{X}\rightarrow \mathcal{X}$ qui v\'{e}rifie les
propri\'{e}t\'{e}s

$(es_{1})$ $\varphi(t,t,x)=x, \ \forall (t,x)\in
\mathbb{R}_{+}\times \mathcal{X}$

$(es_{2})$ $\varphi(t,s,\varphi(s,t_{0},x))=\varphi(t,t_{0},x), \
\forall (t,s),(s,t_{0})\in \mathcal{T}, \ \forall x\in
\mathcal{X}$

\noindent est appel\'{e}e \emph{semiflot d'\'{e}volution} sur
$\mathcal{X}$.
\end{definition}

\begin{definition}\label{cocycle}\rm
Une application $\Phi:\mathcal{T}\times \mathcal{X}\rightarrow
\mathcal{B}(\mathcal{V})$ qui satisfait les conditions suivantes

$(ec_{1})$ $\Phi(t,t,x)=I, \ \forall t\geq0,\ \forall x\in
\mathcal{X}$

$(ec_{2})$ $\Phi(t,t_{0},x)=\Phi(t,s,\varphi(s,t_{0},x))\Phi(s,t_{0},x),
\ \forall (t,s), (s,t_{0})\in \mathcal{T}, \ \forall x\in
\mathcal{X}$

\noindent s'appelle \emph{application cocyclique d'\'{e}volution}
relative au semiflot d'\'{e}volution $\varphi:\mathbb{R}_{+}\times
\mathcal{X}\rightarrow \mathcal{X}$.
\end{definition}

\begin{definition}\label{lses}\rm
Une application $\xi:\mathcal{T}\times \mathcal{Y}\rightarrow
\mathcal{Y}$ d\'{e}finie par
\begin {equation}
\xi(t,s,x,v)=(\varphi(t,s,x),\Phi(t,s,x)v), \ \forall (t,s,x,v)\in
\mathcal{T}\times \mathcal{Y}
\end{equation}
\noindent o\`{u} $\Phi$ est une application cocyclique
d'\'{e}volution relative au semiflot d'\'{e}volution $\varphi$,
est appel\'{e}e \emph{cocycle d'\'{e}volution} sur $\mathcal{Y}$.
\end{definition}

\begin{example}\rm\label{ex_cev}
Soit l'espace m\'{e}trique
$\mathcal{C}(\mathbb{R}_{+},\mathbb{R})=\{f:\mathbb{R}_{+}\rightarrow
\mathbb{R} \ | \ f \ \textrm{continue}\}$ muni de la distance
\begin{equation*}
d(x,y)=\sum_{n=1}^{\infty}\frac{1}{2^{n}}\frac{d_{n}(x,y)}{1+d_{n}(x,y)}
\end{equation*}
o\`{u}
\begin{equation*}
d_{n}(x,y)=\underset{t\in [0,n]}\sup{|x(t)-y(t)|}.
\end{equation*}

Pour tout $n\in \mathbb{N}^{*}$ on consid\`{e}re une fonction
d\'{e}croissante
\[
x_{n}:\mathbb{R}_{+}\rightarrow \left( \frac{1}{2n+1},\frac{1}{2n}\right) ,\
\underset{t\rightarrow \infty}\lim{x_{n}(t)}=\frac{1}{2n+1}.
\]
On note 
\[
x_{n}^{s}(t)=x_{n}(t+s), \ \forall t,s \geq 0
\]
et
soit $\mathcal{X}$ l'adh\'{e}rence dans
$\mathcal{C}(\mathbb{R}_{+},\mathbb{R})$ de l'ensemble
${\{x_{n}^{s},n\in \mathbb{N}^{*},s\in \mathbb{R}_{+}\}}$.
Alors l'application
\begin{equation*}
\varphi:\mathcal{T}\times \mathcal{X}\rightarrow \mathcal{X}, \
\varphi(t,s,x)=x_{t-s}, \ \textrm{o\`{u}} \
x_{t-s}(\tau)=x(t-s+\tau), \ \forall \tau\geq 0
\end{equation*}
est un semiflot d'\'{e}volution sur $\mathcal{X}$.

Si on consid\`{e}re l'espace de Banach $\mathcal{V}=\mathbb{R}^{p}$, $p\geq 1$ muni de la norme 
\[
\left\Vert
(v_{1},...,v_{p})\right\Vert=|v_{1}|+...+|v_{p}|, 
\]
alors
l'application $\Phi:\mathcal{T}\times \mathcal{X}\rightarrow
\mathcal{B}(\mathcal{V})$ d\'{e}finie par
\begin{equation*}
\Phi(t,s,x)(v_{1},...,v_{p})=\left( e^{\alpha_{1}\int_{s}^{t}x(\tau-s)d\tau}v_{1},...,e^{\alpha_{p}\int_{s}^{t}x(\tau-s)d\tau}v_{p}\right) ,
\end{equation*}
o\`{u} on consid\`{e}re $(\alpha_{1},...,\alpha_{p})\in\mathbb{R}^{p}$ fix\'{e},
est une application cocyclique d'\'{e}volution.

Alors $\xi=(\varphi,\Phi)$ est un cocycle d'\'{e}volution sur
$\mathcal{Y}$.
\end{example}

\begin{definition}\label{cevshift}\rm
L'application $\xi_{\gamma}:\mathcal{T}\times \mathcal{Y}\rightarrow
\mathcal{Y}$ d\'{e}finie par
\begin{equation}\label{relcevshift}
\xi_{\gamma}(t,s,x,v)=(\varphi(t,s,x),e^{-\gamma(t-s)}\Phi(t,s,x)v),
 \forall (t,s)\in \mathcal{T},  \forall (x,v)\in \mathcal{Y},
\end{equation}
o\`{u} $\gamma \in \mathbb{R}$, est nomm\'{e}e \emph{cocycle
d'\'{e}volution shift} sur $\mathcal{Y}$.
\end{definition}

\vspace{3mm}

On va noter par la suite
\[
e^{-\gamma(t-s)}\Phi(t,s,x)v=\Phi_{\gamma}(t,s,x)v, \ (t,s)\in
\mathcal{T}, \ (x,v)\in \mathcal{Y}, \ \gamma\in \mathbb{R}.
\]

Dans cet article on \'{e}tudie une classe particuli\`{e}re de
cocycles d'\'{e}volution introduits dans la d\'{e}finition
suivante.

\begin{definition}\label{sm}\rm
Un cocycle d'\'{e}volution $\xi=(\varphi,\Phi)$ est nomm\'{e}
\emph{fortement mesurable} si l'application donn\'{e}e par $t\rightarrow
\left\Vert \Phi(t,s,x)v\right\Vert$
est mesurable, pour tout $(s,x,v)\in \mathbb{R}%
_{+}\times \mathcal{Y}.$
\end{definition}

\vspace{3mm}

On note
\[
\mathcal{F}=\{f:[0,\infty )\rightarrow (0,\infty ) |  \ f \
\textrm{fonction d\'{e}croissante},  \ \underset{t\rightarrow
\infty }{\lim }f(t)=0\}.
\]

Pour approcher la propri\'{e}t\'{e} de l'instabilit\'{e}, on va
proposer la classe suivante des cocycles d'\'{e}volution.

\begin{definition}\label{dc}\rm
Un cocycle d'\'{e}volution $\xi=(\varphi, \Phi)$ est dit avoir 
\emph{d\'{e}croissance exponentielle} s'il existe une fonction
$f\in \mathcal{F}$ telle que
\begin{equation}\label{rel_dc}
\left\Vert \Phi(t+t_{0},t_{0},x)v\right\Vert \geq f(t)\left\Vert v
\right\Vert,  \forall t,t_{0}\geq 0,  \forall (x,v)\in
\mathcal{Y}.
\end{equation}%
\end{definition}

\vspace{3mm}

Nous avons pr\'{e}sent\'{e} une d\'{e}finition assez
g\'{e}n\'{e}rale. Le r\'{e}sultat suivant peut \^{e}tre
consid\'{e}r\'{e} comme une caract\'{e}risation classique pour la
propri\'{e}t\'{e} de d\'{e}croissance exponentielle.

\begin{proposition}\label{caract_dc}
Le cocycle d'\'{e}volution $\xi=(\varphi, \Phi)$ a une
d\'{e}croissance exponentielle si et seulement s'ils existent des
constantes $\widetilde{N}>1$ et $\omega
>0$ telles que
\begin{equation*}
\widetilde{N}\left\Vert \Phi(t+t_{0},t_{0},x)v \right\Vert\geq
e^{-\omega t}\left\Vert v\right\Vert,  \forall t,t_{0}\geq 0,
\forall(x,v)\in \mathcal{Y}.
\end{equation*}%
\end{proposition}

\begin{proof}
\textit{N\'{e}cessit\'{e}}. On consid\`{e}re une fonction $f\in \mathcal{F}$.
Alors il existe $\mu>0$ tels que $f(\mu )<1$.
De plus, il existe $k\in \mathbb{N}$ et $r\in [0,\mu )$ tel que
\[
t=k\mu +r,  \forall t\geq 0.
\]
On obtient les in\'{e}galit\'{e}s
\begin{equation*}
\left\Vert \Phi(t+t_{0},t_{0},x)v\right\Vert \geq f(r)\left\Vert
\Phi(k\mu +t_{0},t_{0},x)v\right\Vert \geq
\end{equation*}
\begin{equation*}
\geq f(r)f(\mu)\left\Vert \Phi((k-1)\mu
+t_{0},t_{0},x)v\right\Vert \geq ...\geq f(r)f(\mu )^{k}\left\Vert
v\right\Vert
\end{equation*}%
qui sont v\'{e}rifi\'{e}es pour tout $t,t_{0}\geq 0$ et tout $(x,v)\in \mathcal{Y}$.

La d\'{e}finition  de la d\'{e}croissance exponentielle pour $\xi$ est alors obtenue si on note
\[
\widetilde{N}=\frac{1}{f(\mu )}>1  \ \textrm{et}  \ \omega
=-\frac{1}{\mu}\ln f(\mu)>0.
\]

\textit{Suffisance}. On consid\`{e}re la fonction donn\'{e}e par $f(t)=e^{-\omega t}$, ce qui ach\`{e}ve la d\'{e}monstration.
\end{proof}

On va pr\'{e}senter plusieurs cat\'{e}gories de propri\'{e}t\'{e}s asymptotiques pour les cocycles d'\'{e}volution. 

\begin{definition}\label{is}\rm
Un cocycle d'\'{e}volution $\xi=(\varphi, \Phi)$ est dit
\emph{instable} s'il existe une application $N:[0,\infty)\rightarrow (1,\infty)$ telle que%
\begin{equation}
N(t)\left\Vert \Phi(t,t_{0},x)v\right\Vert \geq \left\Vert
v\right\Vert ,  \forall (t,t_{0})\in \mathcal{T},  \forall
(x,v)\in \mathcal{Y}.
\end{equation}%
\end{definition}

\begin{definition}\label{eis}\rm
Le cocycle d'\'{e}volution $\xi=(\varphi, \Phi)$ est nomm\'{e}
\emph{exponentiellement instable} s'il existe une fonction
$N:[0,\infty)\rightarrow (1,\infty)$ et une constante $\nu>0$ avec
la propri\'{e}t\'{e}
\begin{equation}
N(t)e^{-\nu (t-s)}\left\Vert \Phi(t,t_{0},x)v\right\Vert\geq
\left\Vert \Phi(s,t_{0},x)v\right\Vert,
\end{equation}%
\noindent pour tout $(t,s),(s,t_{0})\in \mathcal{T}$ et tout
$(x,v)\in \mathcal{Y}$.
\end{definition}

\begin{example}\rm
On consid\`{e}re $\mathcal{X}=\mathbb{R}_{+}$ et
$\mathcal{V}=\mathbb{R}$. Nous d\'{e}finissons l'application
$\Phi:\mathcal{T}\times \mathbb{R}_{+}\rightarrow
\mathcal{B}(\mathbb{R})$ par
\[
\Phi(t,s,x)v=ve^{t-s-2t\sin\frac{\pi}{4}t+2s\sin\frac{\pi}{4}s}, \
(t,s,x,v)\in\mathcal{T}\times\mathbb{R}_{+}\times\mathbb{R},
\]
qui est une application cocyclique d'\'{e}volution.

Alors $\xi=(\varphi,\Phi)$ est un cocycle d'\'{e}volution sur
$\mathcal{Y}=\mathbb{R}_{+}\times\mathbb{R}$ pour tout semiflot
d'\'{e}volution $\varphi:\mathcal{T}\times
\mathbb{R}_{+}\rightarrow \mathbb{R}_{+}$, exponentiellement instable avec
$N(t)=e^{4t}$ et $\nu =3$.
\end{example}

\begin{definition}\label{d_uiis}
Le cocycle d'\'{e}volution $\xi=(\varphi,\Phi)$ est appel\'{e}
\emph{int\'{e}gralement instable} s'il est fortement mesurable et s'il
existe $M:[0,\infty )\rightarrow [1,\infty )$ telle que
la relation suivante est v\'{e}rifi\'{e}e
\begin{equation}\label{ineguiis}
\int_{0}^{t}\left\Vert \Phi(\tau ,t_{0},x)v\right\Vert d\tau \leq
M(t)\left\Vert \Phi(t,t_{0},x)v\right\Vert,
\end{equation}%
pour tout $(t,t_{0})\in \mathcal{T}$ et tout $(x,v)\in \mathcal{Y}.$
\end{definition}

\section{Caract\'{e}risations de l'instabilit\'{e} exponentielle}

\begin{remark}\label{obs2}\rm
Analogue au cas uniforme pr\'{e}sent\'{e} dans
\cite{MeSt_ICNODEA}, un cocycle d'\'{e}volution exponentiellement
instable est instable et poss\`{e}de la propri\'{e}t\'{e} de
d\'{e}croissance exponentielle.
\end{remark}

\vspace{3mm}

D'autres connections qu'on peut \'{e}tablir entre les
propri\'{e}t\'{e}s asymptotiques pr\'{e}sent\'{e}es dans D\'{e}finition \ref{sm}, D\'{e}finition \ref{dc},
D\'{e}finition \ref{is}, D\'{e}finition \ref{eis} et 
D\'{e}finition \ref{d_uiis} sont donn\'{e}es par les propositions
suivantes.

\begin{proposition}
Un cocycle d'\'{e}volution int\'{e}gralement instable
$\xi=(\varphi,\Phi)$ ayant la propri\'{e}t\'{e} de
d\'{e}croissance exponentielle est instable.
\end{proposition}

\begin{proof}
On consid\`{e}re une fonction $f\in \mathcal{F}$ donn\'{e}e comme
dans la D\'{e}finition \ref{dc}.

Ensuite on va noter $K=\int\limits_{0}^{1}f(\tau)d\tau$.

On obtient successivement
\begin{equation*}
K\left \Vert v\right \Vert =\int\limits_{0}^{1}f(\tau)\left \Vert
v\right \Vert d\tau = \int\limits_{t_{0}}^{t_{0}+1}f(u-t_{0})\left
\Vert v\right \Vert du \leq
\end{equation*}
\begin{equation*}
\leq \int\limits_{0}^{t}\left \Vert \Phi(u,t_{0},x)v\right \Vert
du \leq M(t)\left \Vert \Phi(t,t_{0},x)v\right \Vert
\end{equation*}
pour tout $t\geq t_{0}+1 >t_{0}\geq 0$ et tout $(x,v)\in
\mathcal{Y}$, l'existence de la fonction $M:[0,\infty )\rightarrow
[1,\infty )$ \'{e}tant assur\'{e}e par la D\'{e}finition
\ref{d_uiis}.

Maintenant, si $t\in [t_{0},t_{0}+1)$, on a pour tout $(x,v)\in
\mathcal{Y}$
\[
\left \Vert \Phi(t,t_{0},x)v\right \Vert \geq f(t-t_{0})\left
\Vert v\right \Vert \geq f(1)\left \Vert v\right \Vert.
\]

Alors, si on d\'{e}finit
\[
N(t)=[f(1)]^{-1}+K^{-1}M(t),  \forall t\geq 0,
\]
on obtient la conclusion.
\end{proof}

Une caract\'{e}risation int\'{e}ressante pour l'instabilit\'{e}
exponentielle est pr\'{e}sent\'{e}e en faisant appel \`{a} la
D\'{e}finition \ref{cevshift}.

\begin{proposition}
Soit $\xi=(\varphi,\Phi)$ un cocycle d'\'{e}volution fortement
mesurable ayant d\'{e}croissance exponentielle. Alors $\xi$ est
exponentiellement instable si et seulement s'il existe une
constante $\alpha
>0$ telle que le cocycle d'\'{e}volution shift $\xi_{\alpha}=(\varphi,\Phi_{\alpha})$ soit int\'{e}gralement instable.
\end{proposition}

\begin{proof}
\emph{N\'{e}cessit\'{e}.} On d\'{e}finit
\[
\alpha=\frac{\nu}{2}>0,
\]
o\`{u} l'existence de $\nu$ est assur\'{e}e par l'hypoth\`{e}se et
la D\'{e}finition \ref{eis}, ce qui nous permet d'obtenir les
relations suivantes
\begin{equation*}
\int_{0}^{t}\left\Vert \Phi_{\alpha}(s,t_{0},x)v\right\Vert ds
=\int_{0}^{t}e^{-\alpha(s-t_{0})}\left\Vert
\Phi(s,t_{0},x)v\right\Vert ds \leq
\end{equation*}
\begin{equation*}
\leq N(t)\int_{t_{0}}^{t}e^{-\alpha(s-t_{0})}\left\Vert
\Phi(t,t_{0},x)v\right\Vert e^{-\nu (t-s)}ds\leq
\alpha^{-1}N(t)\left\Vert \Phi_{\alpha}(t,t_{0},x)v\right\Vert,
\end{equation*}
pour tout $(t,t_{0})\in \mathcal{T}$ et tout $(x,v)\in
\mathcal{Y}$. 

Donc, l'instabilit\'{e} int\'{e}grale pour le cocycle
d'\'{e}volution shift est prouv\'{e}e.

\emph{Suffisance.} On note 
\[
K=\int\limits_{0}^{1}e^{-\alpha%
u}f(u)du,
\]
o\`{u} la fonction $f\in \mathcal{F}$ satisfait la
relation (\ref{rel_dc}).

On obtient successivement
\begin{equation*}
K\left\Vert v\right\Vert=\int_{t_{0}}^{t_{0}+1}e^{-\alpha
(\tau-t_{0})} f(\tau -t_{0})\left\Vert
\Phi(t_{0},t_{0},x)v\right\Vert d\tau \leq
\end{equation*}%
\begin{equation*}
\leq \int_{t_{0}}^{t_{0}+1}e^{-\alpha (\tau-t_{0})}\left\Vert
\Phi(\tau,t_{0},x)v\right\Vert d\tau \leq M(t)\left\Vert
\Phi_{\alpha}(t,t_{0},x)v\right\Vert=
\end{equation*}
\[
=M(t)e^{-\alpha (t-t_{0})}\left\Vert \Phi (t,t_{0},x)v\right\Vert
\]
pour tout $(t,t_{0})\in \mathcal{T}$ et tout $(x,v)\in
\mathcal{Y}$. 

Alors $\xi$ est exponentiellement instable.
\end{proof}

\section{G\'{e}n\'{e}ralisations du th\'{e}or\`{e}me de Datko dans le cas de l'instabilit\'{e} exponentielle}

Dans cette partie nous allons d\'{e}montrer deux
caract\'{e}risations pour la propri\'{e}t\'{e} d'instabilit\'{e}
exponentielle dans le cas non uniforme. Le premier th\'{e}or\`{e}me implique
l'instabilit\'{e} ainsi que de l'instabilit\'{e} int\'{e}grale.

\begin{theorem}\label{caract_uiis_1}
Un cocycle d'\'{e}volution fortement continu $\xi$ est
exponentiellement instable si et seulement s'il est instable et
int\'{e}gralement instable.
\end{theorem}

\begin{proof}
\textit{N\'{e}cessit\'{e}}. Si $\xi=(\varphi,\Phi)$ est
exponentiellement instable, alors, par la Remarque \ref{obs2},
$\xi$ est aussi instable.

Il existe $N:(0,\infty)\rightarrow (0,\infty)$ et $\nu >0$ tels
que
\begin{equation*}
\int_{0}^{t}\left\Vert \Phi(\tau,t_{0},x)v\right\Vert d\tau\leq
N(t)\left\Vert \Phi(t,t_{0},x)v\right\Vert\int_{0}^{t}e^{-\nu
(t-\tau)} d\tau \leq
\end{equation*}
\begin{equation*}
\leq \nu^{-1}N(t)\left\Vert \Phi(t,t_{0},x)v\right\Vert
\end{equation*}
pour tout $t\geq t_{0}\geq 0$ et tout $(x,v)\in \mathcal{Y}$, et
alors la propri\'{e}t\'{e} de l'instabilit\'{e} int\'{e}grale est
prouv\'{e}e.

\textit{Suffisance}. L'instabilit\'{e} de $\xi$ implique la
d\'{e}croissance exponentielle et ensuite, l'existence d'une
fonction $f\in \mathcal{F}$ telle que, pour tout $(x,v)\in
\mathcal{Y}$, on peut \'{e}crire
\[
\left\Vert \Phi(s,t_{0},x)v\right\Vert\leq
\frac{1}{f(u-s)}\left\Vert \Phi(u,t_{0},x)v\right\Vert,  \forall
t\geq u\geq s\geq t_{0}\geq 0.
\]
D'apr\`{e}s les propri\'{e}t\'{e}s de la fonction $f$ on obtient
\[
(t-t_{0})\left\Vert v\right\Vert \leq
\frac{1}{f(t)}\int_{t_{0}}^{t}\left\Vert
\Phi(u,t_{0},x)v\right\Vert du\leq \widetilde{M}(t)\left\Vert
\Phi(t,t_{0},x)v\right\Vert,
\]
pour tout $(t,t_{0})\in\mathcal{T}$ et tout $(x,v)\in\mathcal{Y}$, o\`{u} on a not\'{e}
\[
\widetilde{M}(t)=\frac{M(t)}{f(t)},  t\geq 0,
\]
l'existence de la fonction $M:[0,\infty)\rightarrow [1,\infty)$
\'{e}tant assur\'{e}e par l'hypoth\`{e}se et par la D\'{e}finition
\ref{d_uiis}. 

Alors, on a prouv\'{e} l'instabilit\'{e} exponentielle de $\xi$.
\end{proof}

Une autre caract\'{e}risation pour la propri\'{e}t\'{e} de l'instabilit\'{e} exponentielle d'un cocycle d'\'{e}volution
est obtenue \`{a} l'aide de la decroissance exponentielle et de l'instabilit\'{e} int\'{e}grale. 

\begin{theorem}\label{caract_uiis_2}
Un cocycle d'\'{e}volution fortement mesurable $\xi$ est
exponentiellement instable si et seulement s'il est
exponentiellement d\'{e}croissant et int\'{e}gralement instable.
\end{theorem}

\begin{proof}
Nous allons d\'{e}montrer premi\`{e}rement qu'un cocycle d'\'{e}volution ayant une d\'{e}croissance exponentielle est instable s'il
existe une application $\widetilde{M}:[0,\infty)\rightarrow
(1,\infty)$ telle que
\begin{equation*}
\widetilde{M}(t)\left\Vert \Phi(t,t_{0},x)v\right\Vert \geq
\left\Vert \Phi(s,t_{0},x)v\right\Vert,  \forall t\geq s+1>s\geq
t_{0}\geq 0,  \forall (x,v)\in Y.
\end{equation*}%

On va consid\'{e}rer une fonction $g$ donn\'{e}e comme dans la
D\'{e}finition \ref{dc}. Alors il existe $\lambda
>1$ tel que $g(\lambda )<1$. Soit $s\geq 0$. Pour tout $t\in [s,s+1)$, par le m\^{e}me
r\'{e}sultat, on obtient les in\'{e}galit\'{e}s suivantes
\begin{equation*}
\left\Vert \Phi(t,t_{0},x)v\right\Vert\geq g(t-s)\left\Vert
\Phi(s,t_{0},x)v\right\Vert \geq g(\lambda )\left\Vert
\Phi(s,t_{0},x)v\right\Vert.
\end{equation*}
Alors, si on note
\[
N(t)=g(\lambda )^{-1}+ \widetilde{M}(t),  t\geq 0,
\]
la propri\'{e}t\'{e} d'instabilit\'{e} de $\xi$ vient d'\^{e}tre
d\'{e}montr\'{e}e.

Soit une fonction $f\in \mathcal{F}$ donn\'{e}e comme dans la
D\'{e}finition \ref{dc}.

On obtient successivement pour tout $t\geq t_{0}+1>t_{0}\geq 0$ et
tout $(x,v)\in \mathcal{Y}$
\begin{equation*}
\left\Vert \Phi(t_{0},t_{0},x)v\right\Vert \int_{0}^{1}f(\tau )
d\tau =\int_{t_{0}}^{t_{0}+1}f(u-t_{0})\left\Vert
\Phi(t_{0},t_{0},x)v\right\Vert du\leq
\end{equation*}%
\begin{equation*}
\leq \int_{t_{0}}^{t_{0}+1}\left\Vert \Phi(u,t_{0},x)v\right\Vert
du\leq \int_{t_{0}}^{t}\left\Vert \Phi(u,t_{0},x)v\right\Vert
du\leq M(t)\left\Vert \Phi(t,t_{0},x)v\right\Vert .
\end{equation*}%

Alors, par le Th\'{e}or\`{e}me \ref{caract_uiis_1}, l'instabilit\'{e}
exponentielle de $\xi$ est d\'{e}montr\'{e}e.
\end{proof}

On va conclure l'\'{e}tude de l'instabilit\'{e} d'un cocycle d'\'{e}volution dans le cas
non uniforme par le corollaire suivant.

\begin{corollary}
Pour un cocycle d'\'{e}volution int\'{e}gralement instable, les
propri\'{e}t\'{e}s de d\'{e}croissance exponentielle,
d'instabilit\'{e} et d'instabilit\'{e} exponentielle sont
\'{e}quivalentes.
\end{corollary}

\begin{proof}
Le r\'{e}sultat est obtenu en faisant appel \`{a} la Remarque
\ref{obs2} et aux Th\'{e}or\`{e}mes \ref{caract_uiis_1} et
\ref{caract_uiis_2}.
\end{proof}

\vspace{5mm}

{\footnotesize

\noindent\begin{tabular}[t]{ll}

Codru\c{t}a Stoica \\
Institut de Math\' ematiques de Bordeaux \\
Universit\' e Bordeaux 1  \\
351 Cours de la Lib\'{e}ration \\
33405 Talence Cedex, France  \\
E-mail: \texttt{codruta.stoica@math.u-bordeaux1.fr}
\end{tabular}

}

\end{document}